# Solution Bounds, Stability and Estimation of Trapping/Stability Regions of Some Nonlinear Time-Varying Systems *

Mark A. Pinsky [1] and Steve Koblik [2]

[1] *M.A. Pinsky, Department of Mathematics & Statistics, University of Nevada, Reno, Reno NV 89557, USA*

[2] *Steve Koblik, 8110 Birchfield Dr. Indianapolis, IN 46268, USA*

Correspondence should be addressed to Mark Pinsky; pinsky@unr.edu

**Abstract.** Estimation of solution norms and stability for time-dependent nonlinear systems is ubiquitous in numerous engineering, natural science and control problems. Yet, practically valuable results are rare in this area. This paper develops a novel approach, which bounds the solution norms, derives the corresponding stability criteria, and estimates the trapping/stability regions for some nonautonomous and nonlinear systems, which arise in various application domains. Our inferences rest on deriving a scalar differential inequality for the norms of solutions to the initial systems. Utility of the Lipschitz inequality linearizes the associated auxiliary differential equation and yields both the upper bounds for the norms of solutions and the relevant stability criteria. To refine these inferences, we introduce a nonlinear extension of the Lipschitz inequality, which improves the developed bounds and allows estimation of the stability basins and trapping regions for the corresponding systems. Finally, we confirm the theoretical results in representative simulations.

## 1. Introduction

This paper derives a scalar differential inequality and corresponding first order nonlinear auxiliary equation that bounds in norm the solutions of nonautonomous nonlinear system,

$$\dot{x} = A(t)x + f(t,x) + F(t), \quad t \in [t_0, \infty), \quad x \in \mathbb{R}^n, \quad f(t,0) = 0 \qquad (1)$$

where functions, $f := [t_0, \infty) \times \mathbb{R}^n \to \mathbb{R}^n$, $F := [t_0, \infty) \to \mathbb{R}^n$ and matrix, $A(t) \in \mathbb{R}^{n \times n}$, $\forall t \geq t_0$ are continuous, $t_0 \in \mathrm{T} := [\varsigma, \infty)$, $\varsigma \in \mathbb{R}$, $F(t) = \hat{F}\phi(t)$, $\|\phi(t)\| = 1$ and a scalar $\hat{F} \geq 0$. Note also that throughout this paper symbol, $\|\cdot\|$ stands for 2-norm unless it is indicated otherwise. To simplify notation, we will write that, $x(t, t_0, x_0) \equiv x(t, x_0)$, where $x(t, t_0, x_0)$ be a solution to the initial value problem for (1), i.e., $x(t_0, t_0, x_0) = x_0$. We assume below that, $x(t, x_0)$ is uniquely defined for $\forall t \geq t_0 \in \mathrm{T}$ and $x_0 \in \Gamma_{x_0} \subset \mathbb{R}^n$, where $\Gamma_{x_0}$ is a neighborhood of $\mathbb{R}^n$ containing, $x_0 \equiv 0$. Note that the pertained conditions can be found, e.g., in [1] and [2].

We also examine the solutions to homogeneous counterpart to (1),

$$\dot{x} = A(t)x + f(t,x) \qquad (2)$$

Development of efficient stability criteria for the trivial solution to (2) is essential in numerous applied and control problems. For instance, these criteria enable the design and analysis of performance of robust controllers and observers [3].

.

There are two main approaches to this problem: the Lyapunov functions method, see for instance [2] and [3], and the first approximation methodology, see e.g. reviews [4] and [5] as well as [6] and [7] for additional references and historical perspectives. The former approach is widespread in control literature, see [2], [8]- [17] and additional references therein. However, adequate Lyapunov functions are rare for time dependent and nonlinear systems.

The latter approach delivers sufficient stability criterion under the following conditions, see [6] and [7]. The first is the Lipschitz condition,

$$\|f(x,t)\| \leq l(t)\|x\|, \quad \forall x \in \Omega_1 \subset \mathbb{R}^n, \forall t \geq t_0 \in T \tag{3}$$

where $\Omega_1$ is a bounded subset in $\mathbb{R}^n$ containing, $x \equiv 0$ and the function, $l(t) > 0$ is continuous, and, $l(t) \leq \hat{l} = const, \quad \forall t \geq t_0 \in T$. The second condition,

$$\|W(t,t_0)\| \leq N e^{-\lambda(t-t_0)}, \quad \forall t \geq t_0 \in T, \ 0 < N, \lambda = const \tag{4}$$

bounds the growth rate of the transition matrix, $W(t,t_0) = W(t)W^{-1}(t_0)$, where $W(t)$ is the fundamental solution matrix of the linearized equation (2),

$$\dot{x} = A(t)x \tag{5}$$

Inequality (4) comprises necessary and sufficient conditions for asymptotic/exponential stability of (5), e.g. [3] and [6]. Consequently, it was shown that the trivial solution to (2) is exponentially stable if (3), (4), and the following condition,

$$N\hat{l} - \lambda < 0 \tag{6}$$

hold [6], [7]. A somewhat more flexible condition on the growth of $W(t,t_0)$ was introduced in [5], see also [4]

$$\|W(t,t_0)\| \leq N \exp \int_{t_0}^{t} \Upsilon(s)ds, \ \forall t > t_0 > 0 \tag{7}$$

where $\Upsilon(t)$ is an integrable function. Clearly, (7) reduces to (4) for $\Upsilon(t) = const$. In turn, (3) and (7) provide asymptotic stability of the trivial solution to (2) if [4],

$$\lim_{t_0 < t \to \infty} \sup \left(1/(t-t_0)\right) \int_{t_0}^{t} \Upsilon(s)ds + N\hat{l} < 0 \tag{8}$$

While the existence of (4) is acknowledged under some broad conditions [6], to our knowledge, there were no attempts to adequately define function, $\Upsilon(t)$ in (7) and to apply either criterion to stability analysis of practically relevant systems. Furthermore, it was shown, e.g. in [18], that the time-histories of different estimates of the Euclidian norms for the second order fundamental matrix, i.e. $\|W(t)\| = \|\exp At\|$, $A = const$, can diverge from each other and the exact values of $\|\exp At\|$. This raises concern of the practical value of the listed above sufficient stability criteria. Furthermore, in section 3, we show that (4) and (7) can be viewed as conservative versions of the estimate of the norm of transition matrix that follows from our approach.

An attempt to escape the utility of prior bounds on $\|W(t,t_0)\|$ in stability analysis of (2) was undertaken in [19]. However, authentication of the developed stability conditions for relatively complex systems can present a challenging task for this approach as well.

The problem of estimating the norms of solutions to (1) subject to (3) and (4) was reviewed in [6] and [7].

The problem of estimating the states of linear and nonlinear systems was considered in [20] – [24].

The contribution of this paper is two-fold. First – is methodological. This paper derives a novel scalar differential inequality for the norms of solutions to a practically important class of systems governing by the equations (1) or (2), which collapses the dimension of the original estimation problem to one. Due to the comparison principle [3] solutions to this inequality are bounded from above by the solutions to the auxiliary scalar first order linear or nonlinear equations with variable coefficients, which are devised and analyzed in this paper. The linear auxiliary equation is obtained via application of Lipschitz condition, whereas the nonlinear auxiliary equation is devised through application of a nonlinear version of Lipschitz condition, which is also derived in this paper. The second contribution is in application of the conceived methodology to various local and nonlocal estimation and stability problems. This includes utility of the linear auxiliary equation in derivations of relaxed and more general local boundedness and stability criteria as well as application of the nonlinear auxiliary equation to estimation of solution bounds and trapping/stability regions of solutions to the original systems. Our approach bypasses utility of Lyapunov functions method. The conceived approach enhances stability criteria (i.e. (6) and (8)) that are devised in the context of Lyapunov first approximation methodology and develops novel stability and boundedness criteria.

Our inferences are validated in simulations of the Van der- Pol like model, which includes a time dependent linear block and oscillatory external force.

This paper is organized as follows. The next section derives the pivotal differential inequality and pertained auxiliary equation. The subsequent section linearizes the auxiliary equation via utility of Lipschitz inequality and develops the corresponding solution bounds and stability criteria. Section 4 introduces a nonlinear extension of the Lipschitz inequality and develops its various applications, section 5 presents the simulation results and section 6 concludes this study.

## 2. Differential Inequality for Solution Norms

This section derives the pivoting scalar differential inequality for 2-norms of solutions to (1) or (2), which is analyzed subsequently in this paper. Note that an attempt to derive directly from (1) a scalar differential equation governing evolution of only $\|x(t, x_0)\|$ fails. In fact, if (1) is written on spherical coordinates, then the equation for the radius-vector, $r(t, x_0) = \|x(t, x_0)\|$ includes also the angel variables that cannot be discounted in general.

Instead, for the broad class of nonlinear systems we derive below the initial value problem including a scalar differential *inequality* for $\|x(t, x_0)\|$ and the matching initial condition for this function. Using the comparison principle [3], we bound from above a solution to this problem by the matching solution to the associated initial value problem for the auxiliary scalar differential equation. Finally, this last solution bounds in norm from above the solution to (1) with consistent initial value. This allows to collapse dimension and drastically simplify the problem of estimating time-histories of the actual norms of solutions to equations (1) or (2).

In fact, the application of variation of parameters lets us derive from (1) the following equation, e.g. [3],

$$x(t, x_0) = W(t)W^{-1}(t_0)x_0 + W(t)\int_{t_0}^{t} W^{-1}(\tau)\left(f(\tau, x(\tau, x_0)) + F(\tau)\right)d\tau, \ \forall t \geq t_0 \in \mathrm{T}$$

where $W(t)$ is frequently normalized to satisfy the condition, $W(t_0) = I$, where $I$ is the identity matrix. In section 5 we present normalization of $W(t)$, which is more natural for our studies and, hence, used consequently in our simulations. Presently, we only assume that $\|W(t_0)\| = 1$. The last equation leads to the following inequality,

$$\|x(t, x_0)\| \leq \|W(t)\| \|W^{-1}(t_0)x_0\| + \|W(t)\| \int_{t_0}^{t} \|W^{-1}(\tau)\| \|f(\tau, x(\tau, x_0)) + F(\tau)\| d\tau, \ \forall t \geq t_0 \in \mathrm{T} \tag{9}$$

Next, we attempt to match the solutions to (9) with the solution to the initial value problem to a scalar inequality that can be written as follows,

$$D^+ X_1 \leq p(t) X_1 + k(t) \|f(t, x(t, x_0)) + F(t)\|, \ \forall t \geq t_0 \in \mathrm{T}$$
$$X_1(t_0) = X_0 \tag{10}$$

where, $X_1(t, X_0) := [t_0, \infty) \times [0, \infty) \to [0, \infty)$ is a continuous function, which is used to mimic $\|x(t, x_0)\|$ in (9), $D^+$ is Dini's upper right-hand derivative in $t$ [3], functions, $p(t)$ and $k(t)$ are from $[t_0, \infty)$ to $\mathbb{R}$ and $[t_0, \infty)$ to $[1, \infty)$, respectively, $X_0$ is a nonnegative scalar and function, $x(t, x_0)$ is a solution to (9). Note that functions, $p(t)$ and $k(t)$, and the initial value, $X_0$ are uniquely defined below via matching the solutions to (10) and the right hand-side of (9).

Due to comparison principle [3, pp.102-104], solutions of inequality (10) are bounded from above by solutions of the matching differential equation,

$$\dot{X}_2 = p(t) X_2 + k(t) \|f(\tau, x(\tau, x_0)) + F(t)\|$$
$$X_2(t_0) = X_0 \tag{11}$$

Hence, $X_1(t, X_0) \leq X_2(t, X_0)$, $\forall t \geq t_0$.

Then, the application of variation of parameters to (11) yields,

$$X_1(t, X_0) \leq e^{\zeta(t)} \left( X_0 + \int_{t_0}^{t} e^{-\zeta(\tau)} k(\tau) \psi(\tau, x(\tau, x_0)) d\tau \right) \tag{12}$$

where, $\zeta(t) = \int_{t_0}^{t} p(s) ds$ and, $\psi(\tau, x(\tau, x_0)) = \| f(\tau, x(\tau, x_0)) + F(\tau) \|$.

Consequently, we determine $p(t)$, $k(t)$ and $X_0$ by matching the right hand-sides of (9) and (12). Comparison of the first additions in the right hand-sides of (9) and (12), i.e., $\|W(t)\| \|W^{-1}(t_0) x_0\|$ and $e^{\zeta(t)} X_0$, returns,

$$\|W(t)\| = \exp\left( \int_{t_0}^{t} p(s) ds \right) \tag{13}$$

$$\|W^{-1}(t_0) x_0\| = X_0 \tag{14}$$

Next, from (13), $e^{-\zeta(t)} = 1/\|W(t)\|$. Equating the last additions in the right hand-side of (9) and (12) and multiplying and dividing the former function by $\|W(\tau)\|$, returns,

$$\int_{t_0}^{t} e^{-\zeta(\tau)} k(\tau) \psi(\tau, x(\tau, x_0)) d\tau = \int_{t_0}^{t} \left( \|W(\tau)\| \|W^{-1}(\tau)\| \right) \psi(\tau, x(\tau, x_0)) / \|W(\tau)\| d\tau$$

The last relation yields that, $k(t) = \|W(t)\| \|W^{-1}(t)\| = \sigma_{\max}(W) / \sigma_{\min}(W)$ is the running condition number of $W(t)$, and $\sigma_{\max}(t)$ and $\sigma_{\min}(t)$ are maximal and minimal singular values of $W(t)$. Note that $\|W^{-1}(t)\| < \infty$, $\forall t \geq t_0$ and $\sigma_{\min}(t) > 0$, $t \geq t_0$ since $W(t)$ is nonsingular matrix and $k(t)$ is continuous since both $\sigma_{\max}$ and $\sigma_{\min}$ are continuous functions due to our initial condition on matrix $A(t)$.

We will assume throughout this paper that, $\sigma_{\max}(t)$ is unique $\forall t \geq t_0 \in \mathrm{T}$. It follows from (13) that,

$$p(t) = d \left( \ln \|W(t)\| \right) / dt$$

is continuous since $\sigma_{\max}$ is unique.

Hence, our definition of $p(t)$, $k(t)$ and $X_0$ implies that solutions to (9) and (10) corresponding to the same $x_0$ are equal each other and that $\|x(t, x_0)\| \leq X_2(t, X_0)$, $\forall X_0 = \|W^{-1}(t_0) x_0\|$, $\forall t \geq t_0$, where $x(t, x_0)$ is a solution to (1) or (2).

Next, multiplication of (13) by $\|W^{-1}(t_0)\|$ and utility of a standard norm's inequality let us rewrite (13) in the following form,

$$\|W(t) W^{-1}(t_0)\| \leq k(t_0) \exp\left( \int_{t_0}^{t} p(s) ds \right) \tag{15}$$

where we use that, $k(t_0) = \|W^{-1}(t_0)\|$ since, $\|W(t_0)\| = 1$. Next, it follows from (9) that $\|x_0\| \leq X_0$ and from (14) that $\|x_0\| \geq X_0 / k(t_0)$. Hence,

$$X_0 / k(t_0) \leq \|x_0\| \leq X_0 \tag{16}$$

Subsequent multiplication of (15) by $\|x_0\|$, returns,

$$\|x(t, x_0)\| \leq k(t_0) \|x_0\| \exp\left( \int_{t_0}^{t} p(s) ds \right)$$

where in the above formula, $x(t, x_0)$ is a solution to linear equation (5).

Next, we examine the relation between the above formulas and the assumptions (4) and (7), which were used prior in stability theory of the corresponding systems. It follows from (15) that scalars $N$ and $\lambda$ in (4) can be interpreted as follows,

$$N = \sup_{t_0 \in T} k(t_0), \quad -\lambda = \sup_{t_0 \in T} \sup_{t > t_0} \left\{ (t - t_0)^{-1} \int_{t_0}^{t} p(s) ds \right\} \tag{17}$$

where we assume that both $N$ and $\lambda < \infty$. Furthermore, in (7) scalar $N$ can be also interpreted as above, whereas, unknown function, $\Upsilon(t)$ can be interpreted as $p(t)$. Hence, our approach allows to explicate unspecified parameters and function in (4) and (7) and disclose the conservative nature of these inequalities. Consequently, we show below that stability conditions (6) and (8) can be interpreted as conservative versions of ones that are derived in section 3 of this paper, see remarks 1 and 2.

To write (11) in the standard form, we introduce a nonlinear extension of Lipschitz inequality,

$$\|f(t, x)\| \leq L(t, \|x\|), \quad \forall x \in \Omega_2 \in \mathbb{R}^n \tag{18}$$

where, $L := [t_0, \infty) \times [0, \infty) \to [0, \infty)$ is a continuous function in both variables and $L(t, 0) = 0$, and $\Omega_2$ is a bounded subset in $\mathbb{R}^n$ containing $x \equiv 0$. Apparently, for polynomial and some other vector-fields, $\Omega_2 \equiv \mathbb{R}^n$, see Section 4.1, which derives (18) for some sets of vector-functions. Clearly, (18) reduces to (3) if $L$ is linear in $\|x\|$.

Afterwards, let us assume that $x_0 \in B_{r_2} \subset \Omega_2, \forall t_0 \in T$, where $B_{r_2}$ is an open ball that is centered at $x \equiv 0$ with conceivably a sufficiently small radius $r_2$. Then, due to continuity in $t$, $x(t, x_0) \in B_{R_2} \subset \Omega_2, \forall t \in [t_0, t_1]$, $t_1 = t_0 + \delta_1$, where $\delta_1 > 0$ be a conceivably small value, $x(t, x_0)$ is a solution to (2) and $B_{R_2}$ is an open ball that is centered at $x \equiv 0$ with radius, $R_2 = R_2(r_2)$. In turn, the solutions to nonhomogeneous equation (1),

$x(t, x_0) \in B_{R_2} \subset \Omega_2, \forall t \in [t_0, t_1]$ if both, $x_0 \in B_{r_2} \subset \Omega_2, \forall t_0 \in T$ and $\hat{F}$ is a sufficiently small value. Clearly, the latter condition is implication of continuity of solutions to (1) in both, $t$ and parameter $\hat{F}$.

Under this last condition, application of (18) to (11) yields the following differential inequality,

$$\dot{X}_3(t, X_0) \leq p(t) X_3(t, X_0) + k(t) \left( L(t, X_3(t, X_0)) + \|F(t)\| \right), \quad \forall x_0 \in B_{r_2}, \ t \in [t_0, t_1], \ t_0 \in T$$
$$X_3(t_0, X_0) = X_0 = \|W^{-1}(t_0) x_0\| \tag{19}$$

where continuous, $X_3(t, X_0) := [t_0, t_1] \times [0, \bar{X}_0] \to [0, \bar{X}_3]$, $\bar{X}_3(t_1, \bar{X}_0) = \sup_{t_0 \leq t \leq t_1} \sup_{X_0 \in [0, \bar{X}_0]} X(t, X_0)$ and $\hat{F}$ is a sufficiently small value.. Clearly, if $\Omega_2 \equiv \mathbb{R}^n$, then the condition, $x_0 \in B_{r_2}$ can be voided and (19) holds $\forall t \geq t_0$ and $\hat{F} \geq 0$.

Note also that (19) is defined $\forall t \geq t_0$ whereas its relation to (11) and, in turn, to (1), is embraced yet for $\forall t \in [t_0, t_1]$.

In turn, due to comparison principle [3], solutions of (19) are bounded by the consistent solutions to the associated differential equation,

$$\dot{X}_4(t, X_0) = p(t) X_4(t, X_0) + k(t) \left( L(t, X_4(t, X_0)) + \|F(t)\| \right), \quad t \in [t_0, t_1], \ X_0 \in [0, \bar{X}_0], \ t_0 \in T$$
$$X_4(t_0, X_0) = X_0 = \|W^{-1}(t_0) x_0\| \tag{20}$$

where $X_4(t, X_0) \geq X_3(t, X_0) \geq X_2(t, X_0) \geq X_1(t, X_0), \forall t \geq t_0$. Let us assume that the initial value problem for equation (20) possesses a unique solution for $\forall t \geq t_0$ and denote that, $\bar{\bar{X}}(\bar{X}_0) = \sup_{t_0 \leq t} \sup_{X_0 \in [0, \bar{X}_0]} X_4(t, X_0)$. Note that, in general, $\bar{\bar{X}}\left(\bar{X}_0\right)$ can be infinity. In section 4, we assume that $B_{\bar{\bar{X}}} \subset \Omega_2$, where $B_{\bar{\bar{X}}}$ be a ball with radius $\bar{\bar{X}}$. This condition implies that solutions to (20) bound in norm from above the solutions to (11) and, in turn, solutions to (1) for $\forall t \geq t_0$.

In the following section we use (3) instead of (18) to linearize (20) in the neighborhood, $\Omega_1 \in \mathbb{R}^n$. This subsequently leads to a scalar, linear and integrable auxiliary equation, which is defined, in general, on a short time-interval. Next we formulate conditions assuring that the solutions of the corresponding equation remain in $\Omega_1$, $\forall t \geq t_0$ and derive some explicit upper bounds for solutions to equation (1) and the corresponding stability criteria for the trivial solution to (2). Finally, we show that stability conditions (6) and (8) can be regarded as conservative counterparts of ones we outlined below.

### 3. Linearization of Auxiliary Equation via Application of Lipschitz Inequality

In analogy with section 2, we assume in this section that $x_0 \in B_{r_1} \subset \Omega_1$, $\forall t_0 \in T$, where $B_{r_1}$ is an open ball that is centered at $x \equiv 0$ with conceivably small radius $r_1$ and $\hat{F}$ is a sufficiently small number. Then, due to continuity in $t$, $x(t, x_0) \in B_{R_1} \subset \Omega_1$, $\forall t \in [t_0, t_1]$, where $x(t, x_0)$ is a solution to (1), $t_1 = t_0 + \delta_1$, $\delta_1 > 0$ is a sufficiently small number and $B_{R_1}$ is an open ball that is centered at $x \equiv 0$ with radius $R_1 = R_1(r_1)$. Next, application of (3) to (11) and utility of the comparison principle [3] let us substitute (20) by a scalar linear equation,

$$\dot{X}(t, X_0) = \left(p(t) + k(t)l(t)\right) X(t, X_0) + k(t) \|F(t)\|, \quad \forall t \in [t_0, t_1], \quad X_0 \in [0, \hat{X}_0], \quad t_0 \in T$$
$$\|X(t_0, X_0)\| = X_0 = \|W^{-1}(t_0) x_0\|$$
(21)

where continuous functions, $X(t, X_0) := [t_0, t_1] \times [0, \hat{X}_0] \to [0, \hat{X}]$, $\hat{X}(t_1, \hat{X}_0) = \sup_{t_0 \leq t \leq t_1} \sup_{X_0 \in [0, \hat{X}_0]} X(t, X_0)$. Let us also define, $X_L(\hat{X}_0) = \sup_{t_0 \leq t} \sup_{X_0 \in [0, \hat{X}_0]} X(t, X_0)$. Clearly, $X_L(\hat{X}_0)$ is infinity if homogeneous counterpart to (21) is unstable. Some inferences with utility of this quantity are drawn in Remark 4.

Note that due to our previous assumptions, $p(t)$, $k(t)$ and $\|F(t)\|$ are continuous functions, which implies that (21) possesses a unique solution for $\forall t \geq t_0$. However, this solution may not bound in norm the corresponding solution to (1) if it exit $\Omega_1$ for some $t > t_1$. These inferences are presented in

**Theorem 1**. Assume that $x_0 \in B_{r_1} \subset \Omega_1$, $\forall t_0 \in T$, $r_1$, $\hat{F}$ and $\delta_1 > 0$ are conceivably sufficiently small, $t_1 = t_0 + \delta_1$, inequality (3) holds and, due to our previous assumptions, $\sigma_{\max}(t)$ is unique and $F(t)$ is continuous. Then,

$$\|x(t, x_0)\| \leq X(t, X_0) = X_h(t, X_0) + X_{nh}(t), \quad \forall t \in [t_0, t_1], \quad t_0 \in T, \quad X_0 = \|W^{-1}(t_0) x_0\|, \quad X_0 \in [0, \hat{X}_0] \quad (22)$$

where, $x(t, x_0)$ is a solution to (1) and, $X(t, X_0)$ is the solution of (21), and

$$X_h(t, X_0) = \|W(t)\| \|W^{-1}(t_0) x_0\| \exp\left(\int_{t_0}^{t} k(s) l(s) ds\right) \quad (23)$$

and

$$X_{nh}(t) = \int_{t_0}^{t} \theta(t, \tau) k(\tau) \|F(\tau)\| d\tau \quad (24)$$

and the transition function, $\theta(t, \tau) = \exp\left(\int_{\tau}^{t} (p(s) + k(s) l(s)) ds\right)$.

**Proof.** Clearly, if $x_0 \in B_{r_1} \subset \Omega_1$ then, due to continuity in $t$, $x(t, x_0) \in B_{R_1} \subset \Omega_1$, $\forall t \in [t_0, t_1]$ and solutions to (1) are bounded in norm from above by the consistent solutions to (21) on the correspondent time interval. The latter linear equation assumes a unique solution, which is defined by (22) - (24). Due to continuity of the underlying functions, the integrals in the last formulas are defined for $\forall t \geq t_0$. □

Hence, the problems of assessing the asymptotic/exponential stability of the trivial solution to (2) or boundedness of solutions to (1) are simplified and comprised in evaluation of the matching problems to the auxiliary linear first order homogeneous/nonhomogeneous equations, and assuring that solutions to (1) or (2), i.e. $x(t, x_0) \in \Omega_1$, $\forall t \geq t_0$. The latter, in turn, can be inferenced under some conditions that are listed below.

Note that the necessary and sufficient conditions for various types of stability of a scalar linear equation is known, e.g. [6], and recently are reviewed in [24], where additional references can be found. Application of these conditions to our first order linear auxiliary equation facilitates development of the matching stability criteria for the trivial solution to nonlinear equation (2). Below we present only some of the most explicit boundedness/stability conditions for equations (1) and (2), which are directly follow from (22) – (24).

Note that subsequent corollaries 1 - 4 assume that conditions of the Theorem 1 and formulas (22)-(24) are embraced and include only the additional conditions that are essential to a specific statement listed below.

**Corollary 1.** Assume that $\|W(t)\| k(t) l(t) + d\|W(t)\|/dt \leq 0$, $\forall t \geq t_* \geq t_0 \in T$, $t_* < t_1$ and the equality in the last formula can be attained only for some isolated values of $t$. Then, the trivial solution to (2) is asymptotically stable.

**Proof.** In fact, since, $\|W^{-1}(t_0)\| < \infty$ and, due to (23), $X_h(t_*, X_0) = O(\|x_0\|) < R_1$ for sufficiently small $\|x_0\|$. Differentiating (23) in $t$ implies that for $\forall t \geq t_*$ either, $dX_h(t, X_0)/dt < 0$ or $dX_h(t, X_0)/dt = 0$ for some isolated values of $t$, which yields that, $X_h(t, X_0) \to 0$ monotonically with $t_* \leq t \to \infty$. Hence, $X_h(t_1, X_0) < X_h(t_*, X_0) < R_1$ and $X_h(t_1, X_0)$ can be made arbitrary small for sufficiently small $\|x_0\|$.

Next, application of (22) yields that, $\|x(t_1, x_0)\| \leq X_h(t_1, X_0) < X_h(t_*, X_0) < R_1$. Then, due to continuity of $x(t, x_0)$ in $t$, (22) can be further applied for $\forall t \in [t_1, t_2]$, $t_2 = t_1 + \delta_2$, where $\delta_2 > 0$ is conceivably a small number. In turn, $\|x(t_2, x_0)\| \leq X_h(t_2, X_0) < X_h(t_1, X_0) < R_1$. Hence, (22) can be further protracted for, $\forall t \in [t_2, t_3]$, $t_3 = t_2 + \delta_3$ with conceivably small $\delta_3 > 0$.

Let us show that repetition of these steps yields that, $t_m \to \infty$ if $m \to \infty$. In fact, if, in contrary, $\lim_{m \to \infty} t_m = \tilde{t} < \infty$, then, $\lim_{m \to \infty} \delta_m = 0$ and $x(\tilde{t}, x_0)$ belongs to the boundary of $\Omega_1$. However, due to continuity of $x(t, x_0)$ in $t$, $\lim_{m \to \infty} \|x(t_m, x_0)\| = \|x(\tilde{t}, x_0)\| \leq \lim_{m \to \infty} X_h(t_m, X_0) = 0$. This contradiction infers that, $t_m \to \infty$ if $m \to \infty$. Consequently, our last relation yields that, $\lim_{m \to \infty} \|x(t_m, x_0)\| = 0$ □

**Corollary 2.** Assume that, $\sup_{t \geq t_*}(p(t) + k(t) l(t)) \leq -\nu$, $\nu > 0$, $\forall t \geq t_* \geq t_0 \in T$, $t_* < t_1$. Then, the trivial solution of (2) is exponentially stable. If, in addition, both $\hat{F}$ and $\|x_0\|$ are sufficiently small and $\sup_{t \geq t_0} k(t) = \hat{k} < \infty$, $\forall t_0 \in T$, then the corresponding solutions to (1) are bounded in norm and, $\limsup_{t \to \infty, t \geq t_*} \|x(t, x_0)\| \leq \hat{F}\hat{k}/\nu$.

**Proof.** Assume firstly that $F = 0$ and let $x(t_*, x_0) = x_*$. Due to Theorem 1, (22) holds, $\forall t \in [t_0, t_1]$. Then application of (13) and replacement of $t_0$ by $t_*$ in (23) yields that,
$\|x(t, x_*)\| \leq \|W^{-1}(t_*)\| \|x_*\| \exp(-\nu(t - t_*))$, $\forall t \in [t_*, t_1]$. Next, we recall that, $\|W^{-1}(t_*)\| < \infty$ and, $\|x_*\| = \gamma \|x_0\|$, $0 < \gamma < \infty$ since $x(t, x_0)$ is continuous in $t$. Thus,
$\|x(t_1, x_*)\| \leq \gamma \|W^{-1}(t_*)\| \|x_0\| \exp(-\nu(t - t_*)) < R_0 < R_1$ if $\|x_0\| < R_0 \exp(\nu(t_1 - t_*))/\|W^{-1}(t_*)\|\gamma$. Due to continuity

of $x(t, x_0)$ in $t$, this last inequality let us to extend application of (22) for $\forall t \in [t_1, t_2]$, $t_2 = t_1 + \delta_2$ with sufficiently small $R_0$ and $\delta_2 > 0$. Hence, $\|x(t_2, x_*)\| \leq R_0 e^{-\nu \delta_2}$. Replication of these steps yields that,

$\|x(t_m, x_*)\| \leq R_0 e^{-\nu \Delta_m}$, $\Delta_m = \sum_{i=2}^{m} \delta_i$. Next, following the steps used in Corollary 1, we infer that $\lim_{m \to \infty} \Delta_m = \infty$, which, in turn, yields that $\|x(t_m, x_*)\| \to 0$ exponentially with $m \to \infty$ and proofs the first part of this statement.

Let us assume next that, $\hat{F} > 0$ and both $\hat{F}$ and $\|x_0\|$ are sufficiently small. Due to the made assumption, $\theta(t, \tau) \leq \exp(-\nu(t - \tau))$, $\forall t \geq \tau \geq t_*$. Thus, $X_{nh}(t) \leq (\hat{F}\hat{k}/\nu)(1 - \exp(-\nu(t - t_*)))$, $\forall t \geq t_*$ and $\sup_{t_* \leq t} X_{nh}(t) \leq \hat{F}\hat{k}/\nu$.

Next, due to first part of this statement, $\lim_{t_* \leq t \to \infty} X_h(t, X_*) = 0$ and $X_h(t, X_*) = O(\|x_0\|)$, $\forall t \geq t_*$. This implies that, $\sup_{t_* \leq t} X(t, X_*) = O(\|x_0\| + \hat{F}) < R_+ < R_1$, where $R_+ > 0$ can be made arbitrary small by the appropriate choice of $\|x_0\|$ and $\hat{F}$.

Henceforth, due to Theorem 1, for sufficiently small $\|x_0\|$ and $\hat{F}$, application of (22) implies that, $\|x(t, x_*)\| \leq X(t, X_*) < R_+ < R_1$, $\forall t \in [t_*, t_1]$. Thus, due to continuity of $x(t, x_*)$ in $t$, application of (22) can be extended for $\forall t \in [t_1, t_2]$, $t_2 = t_1 + \delta_2$ with sufficiently small $\delta_2 > 0$. Consequently, $\|x(t, x_*)\| \leq X(t, X_*) < R_+ < R_1$, $t \in [t_1, t_2]$ and $\sup_{t \in [t_1, t_2]} \|x(t, x_*)\|$ can be made arbitrary small if both $\|x_0\|$ and $\hat{F}$ are sufficiently small. As in Corollary 1, repetition of these steps imply that, $t_m \to \infty$ with $m \to \infty$ and, consequently, $\limsup_{m \to \infty, t_m \geq t_*} \|x(t_m, x_0)\| \leq \limsup_{m \to \infty, t_m \geq t_*} X(t_m, X_*) \leq \hat{F}\hat{k}/\nu$ for sufficiently small $\|x_0\|$ and $\hat{F}$ □

**Remark 1**. Clearly, stability condition (6) can be considered as a conservative version of the conditions of the last statement. The former - can be derived from the latter condition by application of (17) and setting, $t_* = t_0$ and $\nu = 0$. In fact, exponential stability of the trivial solution to (2) is assured if $\nu \geq 0$, $\forall t_0 \in T$, but somewhat more conservative condition, $\nu > 0$, $\forall t_0 \in T$ implies uniform exponential stability of trivial solution to (2). Additionally, the condition of the above statement is evoked only for $\forall t \geq t_*$. This discard behavior of solutions on the initial time-interval, where they can diverge from the fixed solution – a common thesis in stability theory.

To formulate less conservative stability criteria, we evoke the definitions of the characteristic and Lyapunov exponents, which determine the fate of solutions to (2) or (5) if $t \to \infty$, see e.g. [4] and [21]. The characteristic exponents, $\chi_i(x(t, x_0)) = \lim_{t_0 < t \to \infty} \sup(t - t_0)^{-1} \ln\|x(t, x_0)\|$, $i \leq n$ assess the rate of exponential growth/decay of $\|x(t, x_0)\|$ if $t \to \infty$. For a linear system (5), the Lyapunov exponents are defined as,

$\mu_i(t_0) = \lim_{t_0 < t \to \infty} \sup(t - t_0)^{-1} \ln \sigma_i(W(t))$, $i \leq n$, where $\sigma_i$ are the singular values of $W$.

Let us also recall that for linear systems the maximal characteristic and Lyapunov exponents are matched, see e.g. [4].

Firstly, we notice that

$$\lim_{t_0 < t \to \infty} \sup\left((t - t_0)^{-1} \int_{t_0}^{t} p(s, t_0) ds\right) = \lim_{t_0 < t \to \infty} \sup\left((t - t_0)^{-1} \ln \|W(t)\|\right) = \mu_{\max}(t_0)$$

where $\mu_{\max}$ is the maximal Lyapunov exponent of solutions to (5). Hence, (5) is asymptotically stable if, $\mu_{\max}(t_0) < 0$, $\forall t_0 \in T$, see e.g. [25], p.94.

Next, using (23) we calculate the characteristic exponent of $X_h(t, X_0)$ as follows,

$$\chi(t_0) = \lim_{t_0 < t \to \infty} \sup(t-t_0)^{-1}\left(\int_{t_0}^{t} k(s)l(s)ds + \ln\|W(t)\|\right)$$

Consequently, we set, $\chi = \mu_{\max} + \chi_*$, where, $\chi_* = \lim_{t_0 < t \to \infty} \sup(t-t_0)^{-1}\int_{t_0}^{t} k(s)l(s)ds$ be a nonlinear correction to the maximal Lyapunov exponent of (5). Afterward, we convey,

**Corollary 3**. Assume that, $\chi(t_0) < 0$, $\forall t_0 \in T$. Then the trivial solution to (2) is asymptotically stable.

**Proof**. In fact, it follows from (23) that, $X_h(t, X_0) \leq \|W^{-1}(t_0)\|\|x_0\|\exp\int_{t_0}^{t}(p(s)+k(s)l(s))ds$, $\forall t \geq t_0 \in T$. Due to our assumption,

$$X_h(t, X_0) = \exp\int_{t_0}^{t}(p(s)+k(s)l(s))ds \leq D_\varepsilon \exp[-(\rho-\varepsilon)(t-t_0)], \forall t \geq t_0 \in T, \rho = -\chi(t_0), D_\varepsilon > 0$$

and $\varepsilon > 0$ be an arbitrary small value, see e.g. [25], p.18 and pp.93-94. Next, since $\|W^{-1}(t_0)\| < \infty$, $X_h(t, X_0) = O(\|x_0\|)D_\varepsilon \exp[-(\rho-\varepsilon)(t-t_0)]$, $\forall t \geq t_0 \in T$ and $X_h(t, X_0) \to 0$ exponentially with $t_0 \leq t \to \infty$. As in the first part of Corollary 2, (22) implies that, $\|x(t, x_0)\| \leq X_h(t, X_0) < R_\times < R_1$, $\forall t \in [t_0, t_1]$ and $\|x(t_1, x_0)\|$ can be made arbitrary small for sufficiently small $\|x_0\|$. Due to continuity of $x(t, x_0)$ in $t$, we can extend application of (22) on the adjacent time-interval, $[t_1, t_2]$, $t_2 = t_1 + \delta_2$, $\delta_2 > 0$, which, in turn, infers that $\|x(t_2, x_0)\| \leq X_h(t_1, X_0)\exp[-(\rho-\varepsilon)\delta_2]$, $\forall t \in [t_1, t_2]$. Replication of the arguments used in the first part of Corollary 2, yields that $\|x(t_m, x_*)\| \leq X_h(t_1, X_0)e^{-\nu\Lambda_m}$, $\Lambda_m = \sum_{i=2}^{m}\Delta_i$ and $\lim_{m \to \infty}\Lambda_m = \infty$, which, in turn, yields that, $\|x(t_m, x_0)\| \to 0$ exponentially with $m \to \infty$. □

**Remark 2.** Condition (8) can viewed as conservative version of the above statement if we assume that,

$$\Upsilon(t) = p(t), N = \sup_{t \geq t_0} k(t) < \infty$$

and $\hat{l} = \sup_{t \geq t_0} l(t) < \infty$. In this case, the left hand-side of (8) present a conservative upper bound for the maximal characteristic exponent of (2). Corollary 3 enhances this bound and discloses its underlined logic.

Let us assume next that $F \neq 0$ and $X_{nh}$ is given by (24) and $\chi(t_0) < 0$, $\forall t_0 \in T$. Then, $\theta(t,\tau) \leq D_\varepsilon(t_0)e^{-\rho(t-\tau)}$, where both $\rho$ and $D_\varepsilon(t_0)$ are defined above, see [25], p.101. In turn, we assume below that $D_\varepsilon(t_0) \leq \hat{D}_\varepsilon$, $\forall t_0 \in T$. This comprises the following,

**Corollary 4**. Assume that, $\chi(t_0) < 0$, $\forall t_0 \in T$, $\theta(t,\tau) \leq \hat{D}_\varepsilon e^{-\rho(t-\tau)}$ and both, $\hat{F}$ and $\|x_0\|$ are sufficiently small values, and, $\sup_{t \geq t_0} k(t) = \hat{k} < \infty$, $\forall t_0 \in T$. Then, solutions of (1) are bounded in norm and,

$\limsup_{t \to \infty} \|x(t, x_0)\| \leq \hat{D}_\varepsilon \hat{F}\hat{k}/\rho$, where, $x(t, x_0)$ is a solution to (1), $-\rho = \chi + \varepsilon$, $\rho > 0$ and $\varepsilon > 0$ is an arbitrary small value, and a constant, $\hat{D}_\varepsilon > 0$.

**Proof**. In fact, due to (24),

$$X_{nh}(t) \leq \hat{F}\hat{k}\hat{D}_\varepsilon \int_{t_0}^{t} e^{-\rho(t-\tau)}d\tau = \left(\hat{F}\hat{k}\hat{D}_\varepsilon/\rho\right)\left(1-\exp(-\rho(t-t_0))\right), \forall t \geq t_0 \in T$$

and $\sup_{t_0 \leq t} X_{nh}(t, X_0) \leq \hat{F}\hat{k}\hat{D}_\varepsilon / \rho$. Next, due to both Corollary 3 and analogy with second part of Corollary 2, we infer that, $\sup_{t_0 \leq t} X(t, X_0) = O(\|x_0\| + \hat{F}) < R_* < R_1$, where scalar $R_* > 0$ can be made arbitrary small for sufficiently small $\|x_0\|$ and $\hat{F}$. Henceforth, due to (22), $\|x(t, x_0)\| \leq X(t, X_0) < R_* < R_1$, $\forall t \in [t_0, t_1]$. Thus, due to continuity of $x(t, x_0)$ in $t$, application of (22) can be extended for $\forall t \in [t_1, t_1 + \delta_2]$ with sufficiently small $\delta_2 > 0$. As in second part of Corollary 2, repetition of such extensions yields that,

$\limsup_{m \to \infty, t_m \geq t_0} \|x(t_m, x_0)\| \leq \limsup_{m \to \infty, t_m \geq t_0} X(t_m, X_*) \leq \hat{F}\hat{k}\hat{D} / \rho$, $t_0 \in T$ for sufficiently small $\hat{F}$ and $\|x_0\|$ □

**Remark 3**. We notice that the application of stability criteria, which are developed in [20] for a scalar linear system, to our auxiliary equation (21) might lead to somewhat less conservative stability criteria for nonlinear equation (2). These types of augments of the above statements are left out of this paper.

**Remark 4.** The proofs of the Corollaries 1 – 4 are simplified under rather more conservative condition on solutions to (21), i.e. $B_{X_L} \subset \Omega_1$, where $B_{X_L}$ is the ball with radius $X_L$. In turn, this condition infers that the solutions to (1) with corresponding initial vectors remain in $\Omega_1$, $t \geq t_0$. Note that similar condition facilitates our inferences in section 4.2.

## 4. Nonlinear Extension of Lipschitz Inequality and its Applications

### 4.1 Extended Lipschitz Inequality

Though application of Lipschitz inequality is widespread in stability and control theories, e.g., [3]- [5], [14], and [15], its utility frequently lead to over conservative inferences, which also evoke dependence of the Lipschitz constant upon the size of the pertaining neighborhood, i.e., $l = l(\Omega_1)$. A rigorous assessment of the last relation can present a challenging task, which is avoided frequently. Yet, this can affect the accuracy of the pertained results. Additionally, admission of (3) linearizes (10) and abates representation of intrinsically nonlocal nonlinear phenomena like, e.g. estimation of trapping/stability regions for the corresponding systems.

To temper these problems, we introduce in this paper a nonlinear extension of Lipschitz inequality, i.e., (18). In principle, a relatively conservative form of (18) can be readily derived for various commonly used functions. For instance, (18) converts to a global inequality, i.e. $\Omega_2 \equiv \mathbb{R}^n$ for polynomial vector fields or ones, which can be presented as a vector Taylor- polynomial with globally bounded Lagrange error term. In these cases, (18) can be attained in a more conservative, but polynomial form, e.g., by successive applications of the following inequalities: $\|f\|_2 \leq \|f\|_1$ and $|x_m^k| \leq \|x\|_2^k$, $m = 1, \ldots, n$, where $x_m$ is the $m$-th component of vector $x$ and $|\cdot|$ is the absolute value. For instance, let, $f : \mathbb{R}^2 \to \mathbb{R}^2$, $x = [x_1 \ x_2]^T$ and, $f = \begin{bmatrix} x_1 x_2^3 & x_1^2 \end{bmatrix}^T$, then

$$\|f\|_2 \leq \|f\|_1 \leq |x_1||x_2|^3 + |x_1|^2 \leq \|x\|_2^4 + \|x\|_2^2 .$$

If the error term in the polynomial approximation of $f$ is bounded for $x \in \Omega_2$, then (18) is validated in the same neighborhood. Yet, such nonlinear inequality frequently appears to be less conservative than (3) in extended neighborhoods of $x \equiv 0$ and let to better represent the underlying behavior of nonlinear systems.

### 4.2 Solution Bounds and Estimation of Trapping/Stability Regions

This section bonds the norms of solutions to either (1) or (2) and estimates the trapping/stability regions for these equations, respectively. Frequently used definitions of these sets of initial vectors are presented below for convenience.

Definition 1. A compact set of initial vectors, that includes zero-vector, is called a trapping region for equation (1), if condition, $x_0 \in \aleph_1$ implies that $x(t, x_0) \in \aleph_1$, $\forall t \geq t_0$.

Definition 2. An open set of initial vectors, that includes zero-vector, is called a stability region of the trivial solution to (2) if condition, $x_0 \in \aleph_2$ implies that $\lim_{t \to \infty} x(t, x_0) = 0$.

Utility of extended Lipschitz inequality (18) frequently sharpens the estimates of the norms of solutions and lessens their dependence upon the size of the pertaining neighborhood, $\Omega_2$, but leads to analysis of solution to the initial value problem for a nonlinear scalar equation with variable coefficients, i.e. (20), which has close form solutions only in some special cases, e.g. if (1) or (2) are autonomous system. Still qualitative analysis and numerical simulations of solutions to a scalar equation is significantly simplified and offers compelling inferences on behavior of solutions to multidimensional systems (1) or (2).

Note that the last two terms in the right side of (20) are nonnegative, whereas $p(t)X_4$ can be either positive or negative or switch the sign for certain values of $t$. This will be used in further analysis of (20).

Utility of solutions to a scalar equation (20) for estimation of the trapping/stability regions for multidimensional equations (1) or (2) requires to relate one-dimensional and $n$-dimensional initial data sets for the corresponding equations. For this sake we define a close set of initial vectors to either equation (1) or (2) as follows,

$$x_0 \in \omega(t_0, X^*) \mid \|W^{-1}(t_0)x_0\| \leq X^*, \ X^* \in [0, \bar{X}_0], \ t_0 \in \mathrm{T}$$

where the set, $\omega(t_0, X^*) \subset \mathbb{R}^n$ is bounded by ellipsoid, $\omega_\gamma(t_0, X^*) \mid \|W^{-1}(t_0)x_0\| = X^*$, which is centered at $x_0 \equiv 0$. Let the open set, $\omega_-(t_0, X^*) = \omega(t_0, X^*) - \omega_\gamma(t_0, X^*)$. Due to (16), $\omega(t_0, X^*) \subset B_{X^*}$, where $B_{X^*}$ be a ball with radius $X^*$. Note also that, $B_{X^*} \subset B_{\bar{X}_0} \subset \Omega_2$, where $B_{\bar{X}_0}$ is a ball with radius $\bar{X}_0$. This leads to the following,

**Theorem 2**. Assume that, $X^* \in [0, \bar{X}_0]$, equation (20) possesses a unique and bounded solution for $\forall X_0 \in [0, \bar{X}_0]$ and $B_{\bar{\bar{X}}(\bar{X}_0)} \subset \Omega_2$, where $B_{\bar{\bar{X}}}$ is the ball with radius $\bar{\bar{X}}$, see section 2 for definition of $\bar{\bar{X}}$. Then,

$$\|x(t, x_0)\| \leq X_4(t, X^*), \ \forall x_0 \in \omega(t_0, X^*), \ \forall t \geq t_0 \tag{25}$$

where $x(t, x_0)$ and $X_4(t, X^*)$ are solutions to either equation (1) or (2) and (20), respectively.

**Proof**. Clearly, since $B_{\bar{\bar{X}}} \subset \Omega_2$, inequality, $\|x(t, x_0)\| \leq X_4(t, X_0)$ with $\|W^{-1}(t_0)x_0\| = X_0$ holds for $t \geq t_0$. Assume now that $X_0 \leq X^*$. Then, $X_4(t, X_0) \leq X_4(t, X^*)$, $\forall t \geq t_0$ since, due to uniqueness, solution curves to a scalar equation (20) do not intersect, which implies (25) □

**Remark 5**. It follows from theorem 2 that $X(t, X_0)$ increases in $X_0$, $\forall t \geq t_0$. This simplifies simulation of (20).

Inequality (25) enables numerical estimation of the trapping/stability regions for (20), which, in turn, leads to estimation of the corresponding regions for the systems (1) or (2).

Note that in the subsequent statement we assume without repetition that conditions of Theorem 2 hold and include only the additional conditions related to this statement.

**Corollary 5.** Assume that solutions to (20) are subjected to one of the following conditions:

1. $F = 0$ and $\lim_{t_0 \leq t \to \infty} X_4(t, X^*) = 0$ for some $0 < X^* < \bar{X}_0$. Then, the trivial solution to (2) is asymptotically stable and $\omega(t_0, X^*)$ is enclosed in its stability basin.

2. $F \neq 0$ and $X_4(t, X^*) \leq d < \bar{X}_0$, $\forall X^* \leq d$, $\forall t \geq t_0$. Then, $\|x(t, x_0)\| \leq d$, $\forall x_0 \in \omega(t_0, d)$, $\forall t \geq t_0$, where $x(t, x_0)$ is a solution to (1), i.e. $\omega(t_0, d)$ is enclosed in the trapping region of solutions to (1).

**Proof.** This corollary directly follows from the Theorem 2. In fact, assume firstly that $F = 0$ and $X_0 \leq X^*$. Then, since $X_4(t, X^*) \leq \bar{\bar{X}}$, $\forall X^* \in [0, \bar{X}_0]$, $\forall t \geq t_0$ and $B_{\bar{\bar{X}}(\bar{X}_0)} \subset \Omega_2$, inequality, $\|x(t, x_0)\| \leq X_4(t, X_0) \leq X_4(t, X^*)$ holds for $\forall x_0 \in \omega(t_0, X^*)$, $\forall t \geq t_0$ and, due to our assumption, yields that, $\lim_{t_0 \leq t \to \infty} \|x(t, x_0)\| = 0$, $\forall x_0 \in \omega(t_0, X^*)$.

Next, let $F \neq 0$. Then, $\|x(t, x_0)\| \leq X_4(t, X_0) \leq X_4(t, X^*) \leq d$, $X_0 \leq X^* \leq d$, $\forall x_0 \in \omega(t_0, d)$, $\forall t \geq t_0$ □

Clearly, the best estimates of the trapping/stability regions yield the maximal admissible values of $X^*$. These values can be readily assessed in simulations of a scalar equation (20), especially, since $X_4(t, X_0)$ is an increasing function in $X_0$, $\forall t \geq t_0$.

Below, we formulate two complementary analytical approaches for estimating such values of $X^*$, which also enhances comprehension of the qualitative structure of solutions to (20). Consequently, next two subsections outline the techniques that bound or approximate equation (20) by its autonomous and integrable counterparts. The first - replaces all time dependent coefficients in (20) by their superior bounds. This yields an autonomous and integrable counterpart of (20) with solutions that bound from above the solution of (20).

The second techniques averages time-dependent coefficient in (20), which yields an autonomous equation approximating (20) under certain conditions. Both techniques lead to explicit solution bounds and boundedness/stability criteria as well as allow to estimate the trapping/stability regions for both autonomous and time-dependent systems.

### 4.3 Reduction of Auxiliary Equation to Autonomous Form

Taking superior bounds of all time-dependent functions in the right-side of (20) yields it autonomous and integrable counterpart,

$$d\hat{X}/dt = \hat{p}\hat{X} + \hat{k}\hat{L}(\hat{X}) + \hat{k}\hat{F} = Q(\hat{X}, \hat{F}) \quad (26)$$
$$\hat{X}(0) = X_0 = \|W^{-1}(t_0)x_0\|$$

where

$$\hat{k} = \sup_{t_0 \leq t} k(t), \quad \hat{p} = \sup_{t_0 \leq t} p(t), \quad \hat{L}(\hat{X}) = \sup_{t_0 \leq t} L(t, \hat{X}) \quad (27)$$

where $\hat{L}(\hat{X}) \geq 0$. Let us assume that, $\hat{k}, \hat{p}, \hat{F} < \infty$, $\hat{L}(0) = 0$ and $\hat{L}(\hat{X})$ is a continuous function, and (26) admits a unique solution, $\hat{X} = \hat{X}(t, X_0)$, $X_0 \in [0, \bar{X}_0]$, $\forall t \geq t_0$. Clearly, $\hat{X}(t, X_0) \geq X_4(t, X_0)$, $\forall t \geq t_0$.

Then, under assumptions, $B_{\bar{\bar{X}}} \subset \Omega_2$, we infer that,

$\|x(t, x_0)\| \leq X_4(t, X_0) \leq \hat{X}(t, X_0)$, $\forall X_0 \in [0, \bar{X}_0]$, $X_0 = \|W^{-1}(t_0)x_0\|$, $\forall t \geq t_0$, where $x(t, x_0)$ is a solution to (1).

As is known, the nonnegative roots of algebraic equation,

$$Q(\hat{X}, \hat{F}) = 0 \quad (28)$$

i.e., $\hat{X} = d_i$, $i = 1, ..., N$ split the initial values of solutions to (26) into subsets, which stem solutions with different behavior on long time-intervals. Subsequently, equations,

$$\|W^{-1}(t_0)x_0\| = d_i \quad (29)$$

link these point-wise boundaries to $n-1$-dimensional ellipsoids in the phase space of (1) or (2), i.e., $\omega_\gamma(t_0, d_i)$, which, in turn, can be used for estimation of the trapping/stability regions for the conforming systems if $\omega(t_0, d_i) \subset \Omega_2$. The latter condition can be assured under the additional assumption, i.e. $B_{\hat{X}_+} \subset \Omega_2$, where

$\hat{X}_+(\bar{X}_0) = \sup_{t_0 \leq t} \sup_{X_0 \in [0, \bar{X}_0]} \hat{X}(t, X_0)$ and $B_{\hat{X}_+}$ be the ball with radius $\hat{X}_+$, which is centered at $x \equiv 0$. Note also that, $B_{\bar{\bar{X}}} \subseteq B_{\hat{X}_+}$.

Below, we review the application of this procedure to some characteristic, but relatively simple cases.

Let us recall that all terms in (26), except, $\hat{p}$ are nonnegative scalars, whereas, $\hat{p}$ can be either positive or negative.

Firstly, we assume for simplicity that $\Omega_2 \equiv \mathbb{R}^n$ and $\hat{p} > 0$. Then $Q(\hat{X},\hat{F}) \geq 0$ and, $\hat{X}(t,X_0) \to \infty$ monotonically with $t \to \infty$ if either, $\hat{F} \neq 0$ and $0 \leq X_0 \leq \bar{X}_0$, or $\hat{F} = 0$ and $0 < X_0 < \bar{X}_0$. Yet, in both cases the norms of the corresponding solutions to (1) or (2) can either approach positive infinity, zero, or remain to be bounded.

Assume next that $F = 0$, $\hat{p} < 0$ and (28) has one simple root, $\hat{X} = d$, $0 < d < \bar{X}_0$. Such fixed solution to (26) can be either stable or unstable. This yields the following,

**Theorem 3**. Assume that $B_{\hat{X}_+} \subset \Omega_2$, $F = 0$, $\hat{p} < 0$ and $\hat{X} = d$, $0 < d < \bar{X}_0$ be a unique fixed solution to (26) corresponding to a simple root of (28)
If this solution is unstable, then the trivial solution to (2) is asymptotically stable. Furthermore,
$\|x(t,x_0)\| \leq d$, $\forall x_0 \in \omega(t_0,d)$, $t \geq t_0$ and, $\lim_{t_0 \leq t \to \infty} \|x(t,x_0)\| = 0$, $\forall x_0 \in \omega_-(t_0,d)$, i.e. $\omega_-(t_0,d)$ is enclosed into the stability basin of the trivial solution to (2).
If this solution is stable, then $\|x(t,x_0)\| \leq d$, $\forall x_0 \in \omega(t_0,d)$, i.e. $\omega(t_0,d)$ is enclosed into the trapping region of the trivial solution to (2) and $\lim_{t_0 \leq t \to \infty} \sup \|x(t,x_0)\| \leq d$, $\forall x_0 \mid \|W^{-1}(t_0)x_0\| = X_0 \in [0,\bar{X}_0]$, where $x(t,x_0)$ is a solution to (2).

**Proof**. The proof of this statement immediately follows from (25) and the assessment of behavior of solutions to a scalar and autonomous nonlinear equation (26) in these two cases. In fact, assume that $F = 0$ and $\hat{X} = d$ be an unique unstable fixed solution to (26). Then, $\hat{X} \equiv 0$ is a stable solution to (26) that attracts all solutions of this equation with initial values, $X_0 < d$, which are monotonically approach zero. Next, since $B_{\hat{X}_+} \subset \Omega_2$, inequality, $\|x(t,x_0)\| \leq \hat{X}(t,X_0)$ holds for $\forall X_0 \in [0,\bar{X}_0]$, $X_0 = \|W^{-1}(t_0)x_0\|$, $\forall t \geq t_0$, where $x(t,x_0)$ is a solution to (2). In turn, this implies that, $\|x(t,x_0)\| \leq d$, $\forall x_0 \mid \|W^{-1}(t_0)x_0\| = X_0 \in [0,\bar{X}_0]$, $t \geq t_0$ and that
$\lim_{t_0 \leq t \to \infty} \|x(t,x_0)\| = 0$, $\forall x_0 \in \omega_-(t_0,d)$. This comprises the first part of the above statement.

Assume next that, $\hat{X} = d$ be a unique stable fixed solution to (26). Then, $\hat{X} \equiv 0$ is the unstable solution to this scalar equation, which implies that both, $\hat{X}(t,X_0) \leq d$, $\forall X_0 \leq d$, $\forall t \geq t_0$ and $\lim_{t_0 \leq t \to \infty} \|\hat{X}(t,X_0)\| = d$, $\forall X_0 \in [0,\bar{X}_0]$. Then, since $B_{\hat{X}_+} \in \Omega_2$, application of inequality $\|x(t,x_0)\| \leq \hat{X}(t,X_0)$, $\forall X_0 \in [0,\bar{X}_0]$, $X_0 = \|W^{-1}(t_0)x_0\|$, $\forall t \geq t_0$ infers that $\|x(t,x_0)\| \leq d$, $\forall x_0 \in \omega(t_0,d)$, $t \geq t_0$ and, in turn, that,
$\lim_{t_0 \leq t \to \infty} \sup \|x(t,x_0)\| \leq d$, $\forall x_0 \mid \|W^{-1}(t_0)x_0\| = X_0 \in [0,\bar{X}_0]$ $\square$

Next, we assume that $\hat{F} > 0$, $\hat{p} < 0$ and (28) has two positive simple roots, $d_i, i = 1,2$, which can be either equal or distinct. Let, $d_1 > d_2$, then we can readily show that, $\hat{X} = d_1$ and, $\hat{X} = d_2$ are unstable and stable fixed solutions to (26), respectively. In fact, in this case, $Q(0,\hat{F}) = \hat{F} > 0$. Hence, continuous, $Q(\hat{X},\hat{F}) > 0$, $\forall \hat{X} < d_2$. Since $\hat{X} = d_2$ is simple root of (27) corresponding to the attractive fixed solution to (26) and, consequently, $\hat{X} = d_1$, is the repelling solution to this equation. This comprises,

**Theorem 4.** Let $B_{\hat{X}_+} \in \Omega_2$, $\hat{F} > 0$, $\hat{p} < 0$ and $x(t,x_0)$ be a solution to (1). In addition, assume that one of the following two conditions hold:
(28) has only two simple roots, $d_i, i = 1,2$ corresponding to unstable and stable fixed solutions to (26), respectively, and $0 < d_2 < d_1 < \bar{X}_0$. Then,

$$\|x(t,x_0)\| \le \begin{cases} d_1, & \forall x_0 \in \omega(t_0,d_1) \\ d_2, & \forall x_0 \in \omega(t_0,d_2) \end{cases}, \forall t \ge t_0,$$

$$\lim_{t_0 \le t \to \infty} \sup \|x(t,x_0)\| \le d_2, \forall x_0 \in \omega_-(t_0,d_1)$$

Let $d_1 = d_2 = d$, then, $\|x(t,x_0)\| \le d, \forall x_0 \in \omega(t_0,d), \forall t \ge t_0$.

**Proof**. The proof of this statement immediately follows from (25) and the assessment of behavior of solutions to a scalar equation (26) in the corresponding cases. In fact, as prior, the condition, $B_{\hat{X}_+} \in \Omega_2$ implies that,

$\|x(t,x_0)\| \le \hat{X}(t,X_0), \forall X_0 \in [0,\bar{X}_0], X_0 = \|W^{-1}(t_0)x_0\|, \forall t \ge t_0$, where $x(t,x_0)$ is a solution to (1). Next, assume firstly that (28) has only two simple roots, $d_1 > d_2 > 0$ corresponding to unstable and stable fixed solutions to (26). Then, since solutions to (26) do not intersect due to uniqueness, $\hat{X}(t,X_0) \le d_i, \forall t \ge t_0$, $i=1,2$ if $0 \le X_0 \le d_i$ and,

$\lim_{t_* \le t \to \infty} \|\hat{X}(t,X_0)\| = d_2, \forall X_0 < d_1$. This with (25) yields the first part of the above statement.

Assume next that, $d_1 = d_2 = d$. In this case, $\hat{X}(t,X_0) \le d$, $0 \le X_0 \le d$, $\forall t \ge t_0$, which with (25) implies the second part of this theorem □

Obviously, (28) can admit more than two positive solutions if $\hat{p} < 0$ and, in addition, the corresponding fixed solutions to (26) can bifurcate due to variation of parameters of this equation. Yet, the corresponding analysis can be extended on these more complex cases alike.

### 4.4 Approximation of Auxiliary Equation Using Averaging Technique

For systems with time dependent linear part, $p(t)$ frequently can be regarded as a highly oscillatory function, i.e., $p = p(t/\varepsilon)$, $0 < \varepsilon \ll 1$, where $\varepsilon$ is a characteristic time scale. Let us also assume for simplicity that $k = k(t/\varepsilon)$, $L = L(t/\varepsilon, X)$, $F = F(t/\varepsilon)$, $t_0 = 0$ and introduce fast time, $\tau = t/\varepsilon$, transforming (20) into the form,

$$dX_5/d\tau = \varepsilon\left(p(\tau)X_5 + k(\tau)L(\tau,X_5)\right) + \varepsilon k(\tau)\|F(\tau)\| = \varepsilon q(\tau,X_5,F(\tau)) \tag{30}$$
$$X_5(\tau = 0, X_0) = X_0, X_0 \in [0,\bar{X}_0]$$

where $X_5(\tau,X_0) \equiv X_4(t = \varepsilon\tau, X_0)$. Since equation (30) is a replica of (20), it admits a unique solution, $X_5(\tau,X_0) \le \bar{\bar{X}}, \forall X_0 \in [0,\bar{X}_0], \forall \tau > 0$ due to our prior assumption on (20). Furthermore, if $B_{\bar{\bar{X}}} \subset \Omega_2$, solutions to (30) with $\tau = t/\varepsilon$ bound in norm from above the corresponding solutions to either (1) or (2), respectively.

Next, application of the averaging technique to (30) yields an autonomous approximation to this equation that can be written as follows,

$$d\tilde{X}_5/dt = \tilde{p}\tilde{X}_5 + \tilde{k}\tilde{L}(\tilde{X}_5) + \tilde{k}\tilde{F} = q_0(\tilde{X}_5, \tilde{F})$$
$$\tilde{X}_5(0) = X_0 = \|W^{-1}(t_0)x_0\| \tag{31}$$

where $\tilde{p} = \lim_{T \to \infty} T^{-1}\int_0^T p(s)ds$, $\tilde{k} = \lim_{T \to \infty} T^{-1}\int_0^T k(s)ds$, $\tilde{F} = \lim_{T \to \infty} T^{-1}\int_0^T \|F(s)\|ds$, $\tilde{L}(\tilde{X}_5) = \lim_{T \to \infty} T^{-1}\int_0^T L(s,\tilde{X}_5)ds$. We assume that the first three limits exist, and the last limit exist uniformly in $\tilde{X}_5 \in [0,\bar{\bar{X}}]$.

Sufficient conditions for the closeness of some solutions of the averaged and initial equations on large and infinite time-intervals can be found in [3], [6], [26], [27], and references therein. For instance, for $\tau \in (0,\infty)$ the following conditions imply closeness of some solutions to (30) and (31), see [6, sec. 7.7].

**Proposition**. Let equation (31) possesses a positive fixed solution, $\tilde{X}_5 = d$ and $B_r(d)$ be a ball with radius $r$, which is centered at $\tilde{X}_5 = d$. Assume that function, $q(\tau, X_5, F(\tau))$ in (30) admits the following conditions:

$q(\tau, X_5, F(\tau)) \leq \infty$, $\forall \tau \in \mathbb{R}$, $X_5(\tau) \in B_r$, and the limit, $q_0(\tilde{X}_5, \tilde{F}) = \lim_{T \to \infty} T^{-1} \int_{\tau}^{\tau+T} q(s, X_5, F(\tau)) ds$, exists uniformly, $\forall X_5 \in B_r$ and $\forall \tau \in \mathbb{R}$, and

$\|\partial q(\tau, X_5, F(\tau))/\partial X_5\| \leq g_1 < \infty$, $\|dq_0(\tilde{X}_5, \tilde{F})/d\tilde{X}_5\| \leq g_2 < \infty$, $\forall \tau \in \mathbb{R}$, $X_5, \tilde{X}_5 \in B_r$.

$\partial^k q(\tau, X_5, F(\tau))/\partial X_5^k$ and $\partial^k q_0(\tilde{X}_5, \hat{F})/\partial \tilde{X}_5^k$, $k = 1, 2$, $\forall \tau \in \mathbb{R}$, $\tilde{X}_5 \in B_r$ are continuous functions.

$q_0^k(\tilde{X}_5, \tilde{F}) = \lim_{T \to \infty} T^{-1} \int_{\tau}^{\tau+T} \left( \partial^k q(\tau, \tilde{X}_5, F(\tau))/\partial \tilde{X}_5^k \right) ds$, $k = 1, 2$ are uniformly defined for $\forall \tau \in \mathbb{R}$, $\tilde{X}_5 \in B_r$.

$\dfrac{dq_0}{d\tilde{X}_{5|\tilde{X}_5=d}} \neq 0$.

Then, for sufficiently small $\mu$, there is $\varepsilon_0 > 0$, such that for $0 < \varepsilon < \varepsilon_0$, (30) admits a unique solution, $z = z(\tau)$, which obeys the inequality,

$$\sup_\tau \|z(\tau) - d\| \leq \mu, \quad \forall \tau > 0 \tag{32}$$

and $z(\tau)$ is a stable/unstable solution to (30) if $\tilde{X}_5 = d$ is a stable/unstable solution to (31).

Yet, under assumption of the Proposition positive solutions, $z_i(\tau)$ belong to $\mu_i$ - neighborhoods of the fixed solutions to (31), $\tilde{X}_5 = d_i$ and assume their stability/instability properties. This let us to use thee fixed solutions for estimation of the trapping/stability regions of (30) nearly the same way it was done prior with utility of the fixed solutions to (26). Subsequent inferences of the behavior of solutions to (1) or (2) can be made under additional condition that, $B_S \subset \Omega_2$, where $B_S$ be a ball with radius $S$, which is centered at zero, $S = \bar{\bar{X}} + s$ and $s = \max_i \mu_i$. Consequently, we note that, $B_{\bar{X}} \subset B_S \subset \Omega_2$.

Next, due to (32), equation (29) should be attuned into the following more conservative form,
$$\|W^{-1}(0) x_0\| = d_i - \mu_i$$

Hence, the ellipsoids, $\omega_\gamma(0, d_i - \mu_i)|$, estimate the sets of the corresponding initial vectors for solutions to (1) or (2), which are included in the trapping/stability regions of these equations if $B_S \subset \Omega_2$. This let us to adjust the statements of theorems 3 and 4 as follows,

**Theorem 5.** Let $\tilde{X}_5(\tau, X_0)$, $\forall \tau \in [0, \infty)$, $X_0 \in [0, \bar{X}_0]$ be a unique solution to (31), $B_S \subset \Omega_2$, $\varepsilon \ll 1$ and functions, $q(\tau, X_5, F(\tau))$ and $q_0(\tilde{X}_5, \tilde{F})$ admit the conditions of the Proposition. Assume also that:

1. $F = 0$, $\tilde{p} < 0$ and (31) possess a unique unstable fixed solution, $\tilde{X} = d \in (0, \bar{X}_0)$. Then, the trivial solution to (2) is asymptotically stable. In addition, $\|x(t, x_0)\| \leq d + \mu$, $\forall x_0 \in \omega(0, d - \mu)$, $t \geq 0$, and $\lim_{t \to \infty} \|x(t, x_0)\| = 0$, $\forall x_0 \in \omega_-(0, d - \mu)$.

If, in turn, $\tilde{X} = d$ be an unique stable fixed solution to (31) then, $\|x(t, x_0)\| \leq d + \mu$, $\forall x_0 \in \omega(0, d - \mu)$, $t > 0$ and $\limsup_{t \to \infty} \|x(t, x_0)\| \leq d + \mu$, $\forall x_0 \| \|W^{-1}(t_0) x_0\| = X_0 \in [0, \bar{X}_0]$, where $x(t, x_0)$ is a solution to (2).

$\tilde{F} > 0$, $\tilde{p} < 0$, and (31) possesses unstable and stable simple fixed solutions, $d_i, i = 1, 2$ and, $0 < d_2 < d_1 < \bar{X}_0$. Then,

$$\|x(t, x_0)\| \leq \begin{cases} d_1 + \mu, & \forall x_0 \in \omega(0, d_1 - \mu) \\ d_2 + \mu, & \forall x_0 \in \omega(0, d_2 - \mu) \end{cases}, \forall t \geq 0$$

$$\limsup_{t \to \infty} \|x(t, x_0)\| \leq d_2 + \mu, \forall x_0 \in \omega_-(0, d_1 - \mu)$$

$\tilde{F} > 0$, $\tilde{p} < 0$, and (31) possesses a repeated fixed solution, $d_1 = d_2 = d < \bar{X}_0$. Then, $\|x(t, x_0)\| \leq d$, $\forall x_0 \in \omega(0, d - \mu)$, $\forall t \geq 0$, where $x(t, x_0)$ is a solution to (1).

**Proof**. The proof of this theorem immediately follows from (32), which infers modifications of Theorems 3 and 4 in the considered cases.

We notice that theoretical estimates for admissible values of $\mu$ and $\varepsilon(\mu)$ turn out to be quite conservative [14], but more accurate estimates frequently can be obtained in numerical simulations.

**Remark 6**. Note that Theorem 5 offers the most lucid application of averaging approach to analysis of solutions to equation (20). Yet, application of averaging technique to (20) with two significantly different time scales yields the equation possessing only slow time, see [26] and [27] and more references therein. It was shown in [16] that under some conditions stability of the system averaged over fast time implies stability of the original system with two-time scales. These inferences can be applied to (20) in the corresponding cases. Moreover, after averaging over fast-time, slow-varying coefficients in (20) frequently can be effectively bounded, which allows efficient convergent of (20) to its more conservative but time-invariant and integrable form.

**Remark 7.** Equation (20) turns into the integrable Bernoulli equation [28] if $F = 0$ and $f(t,x)$ obeys Holder's inequality, i.e., $\|f(t,x)\| \leq c(t)\|x\|^{\alpha}$, $c, \alpha > 0, x \in \Omega \subset \mathbb{R}^n$, which streamlines stability analysis and estimation of the solution bounds for such equations.

**Remark 8**. For the sake of completeness, we briefly compare application of Lyapunov functions methodology aided by one-sided Lipschitz condition with the developed above approach. The former methodology can offer less conservative stability and stabilization conditions of an equilibrium of nonlinear systems than its classical counterpart since the Lipschitz constant is larger or equal to its one-sided analog, see [29] – [32] and additional references therein. For some functions, one-sided Lipschitz constant can assume zero or negative values, which, in principle, can significantly reduce conservatize of the underlined methodology. Nonetheless, the choice of the Lyapunov functions also affect the outcomes of such combined approach. However, efficient Lyapunov functions are rarely available for nonautonomous and nonlinear systems. In contrast, our methodology does not rest on utility of Lyapunov functions.

It was shown in [30] and [31] that application of one-sided Lipschitz condition to the Lyapunov aided design of nonlinear observers is simplified under additional so-known quadratic inner-boundedness condition, which is enforced on nonlinear components of the underlined systems. Utility of both conditions decrease efficacy of this approach and involve estimation of three parameters that depend upon the size of the equilibrium's neighborhood. The computational burden of such task quickly increases in higher dimensions. In contrary, our nonlinear version of Lipschitz inequality can be readily devised in higher dimensions for polynomial vector- fields or ones that can be represented by vector-polynomials with locally/globally bounded error terms. Additionally, the extended Lipschitz inequality becomes global for polynomial vector -fields and if the error terms in polynomial approximations of the vector-field is globally bounded.

Furthermore, our approach naturally enables estimation of trapping/stability regions for time-varying nonlinear systems, whereas the methodology based on utility of one-sided Lipschitz condition, as its classical counterpart, has been primarily used in local stability analysis and stabilization problems.

Finally, we contrast application of both methodologies to some standard cases. Assume, for instance, that $f = x^m$, $x \geq 0$, $m \in \mathbb{R}$. Then, $\|f\| = x^m = L(t,\|x\|)$. Hence, in this case (18) is global and sharp inequality. In turn, one-sided Lipschitz inequality returns, $xf = x^{m+1} \leq \rho x^2$, $\rho = const$. The last inequality fails if $m \neq 1$. Additionally, let $m > 1$ and $x \in [0, x_*]$. Then, one-sided Lipschitz inequality yields that, $\rho = x_*^{m-1}$ and, subsequently, $x^m \leq x_*^{m-1} x$, $\forall x \in [0, x_*]$. The last inequality is sharp if $x = x_*$, otherwise it is inferior than (18).

The above inferences can be extended on polynomial vector-fields. Apparently, less conservative one-sided Lipschitz inequality bears some of shortcomings of its classical counterpart. Nonetheless, one-sided Lipschitz inequality can deliver superior estimates in some application domains.

## 5. Simulations

This section initially applies the developed above methodology for estimating the solution norms as well as trapping/stability regions of Van der-Pol- like model with both time-varying linear part and external time-dependent perturbation. The system is written in dimensionless variables as follows,

$$\dot{x} = A(t)x + f(t,x) + F(t)$$
$$x(0) = x_0 \qquad (33)$$

where $x = [x_1\ x_2]^T$, $A = \begin{pmatrix} 0 & 1 \\ -\omega^2 & -\alpha_1 \end{pmatrix}$, $f = \begin{pmatrix} 0 & -\alpha_2 x_2^3 \end{pmatrix}^T$, and $F(t) = [0\ F_2(t)]^T$, $\omega^2(t) = \omega_0^2 + \omega_1(t)$, $\omega_0 = const$, $\omega_1 = a_1 \sin r_1 t + a_2 \sin r_2 t$, $F_2 = a \sin \omega_2 t$; $a$, $\omega_2$, $a_i$, $r_i = const$, $i = 1,2$. In all further simulations we set $\omega_0 = 2$ and $\alpha_1 = 0.2$. The fundamental solution matrix of the linearized system (33), in general, cannot be found in a closed form. Consequently, we gage the running condition number of such matrix in simulations, which show that $k(t)$ oscillates within certain fixed interval about its mean value, see, e.g. Fig. 4 below.

Firstly, we notice that time-histories of $p(t)$, which are ubiquitous in our analysis, are affected by the normalization of $W(t)$. In instance, we tested two different normalizations: 1). $W(0) = I$ and 2) $W(0) = W_0(0)$, where $W_0(t)$ is the fundamental matrix of solutions of the system, $\dot{x} = Ax$ with $\omega_1 = 0$ and $\|W_0(0)\| = 1$. In Fig. 1, black and blue lines plot time-histories of $p(t)$ corresponding to the first and second normalizations of $W(t)$, respectively. The red line plots running time-average of $p(t)$ for the first normalization of $W(t)$, i.e., $\bar{p}(t) = t^{-1} \int_0^t p(s)\,ds$. It follows from the Fig.1 that the first normalization yields widely oscillatory $p(t)$, but the second – yields, $p(t) = \mathrm{Re}(eig(A)) = -\alpha_1/2$, whereas, $\lim_{t \to \infty} \bar{p}(t) = -\alpha_1/2$. Further simulations show that the second normalization also reduces the variability of $p(t)$ if $\sup|\omega_1|$ assumes relatively small or intermediate values. Hence, the second normalization is adopted in our simulations throughout.

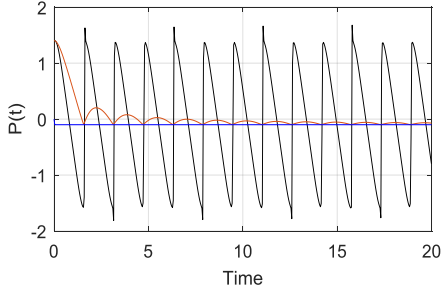

Fig.1. Black and blue lines plot time-histories of $p(t)$ computed for (33) in the following two cases: 1. $W(0) = I$, 2. $W(0) = W_0(0)$, respectively. Red line plots running average of $p(t)$ computed for the first normalization of $W(t)$.

The estimations of the norms of solutions to (33) are shown in Fig.2.1 and Fig.2.2, where red, blue, and black lines plot time-histories of the actual norms of solutions and their two upper bounds, which are comprised either by utility of (3) or nonlinear extension of Lipschitz inequality, i.e. $\left\| \begin{matrix} 0 \\ \alpha_2 x_2^3 \end{matrix} \right\| \le |\alpha_2| \|x\|^3$, $x \in \mathbb{R}^2$, respectively. Hence, for system (33), $\Omega_2 \equiv \mathbb{R}^2$. Yet, the value of Lipschitz constant in (3) depends upon $\sup x_2(t)$ attaining in these simulations. This value in our simulations is estimated using energy integral for the linearized, time invariant and homogeneous model of (33).

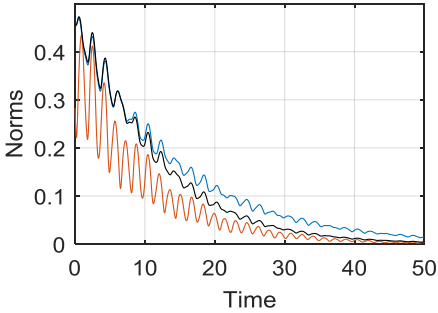
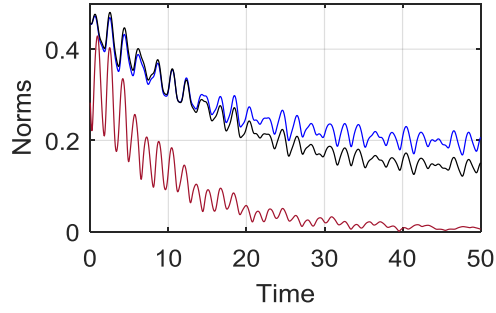

Fig. 2.1    Fig. 2.2

Figure 2. Estimation of solution norms for system (33). Red, blue and black lines plot time-histories of the actual norm of solution to (33) and its upper bounds comprised by utility of Lipschitz or extended Lipschitz inequalities, respectively.

The plots on Fig.2.1 and Fig.2.2 correspond to, $\alpha_2 = 0.1$, $\omega_2 = 2\pi$, $a_1 = a_2 = 0.5$, $r_1 = \pi$, $r_2 = 7$; and in Fig. 2.1, $a = 0$, and in Fig. 2.2, $a = 0.01$.

Clearly, time – histories of the solution bounds comprising the nonlinear extension of Lipschitz inequality outperform the ones utilizing (3) everywhere except a small initial time interval, where the latter is somewhat more accurate than former. Both bounds provide superior accuracy on the initial time intervals, which, however, decreases when time elapses. Application of nonlinear extension of Lipschitz inequality delivers tolerable accuracy on extended time intervals for the homogeneous system. Yet, the estimation accuracy declines for the nonhomogeneous system.

We notice that the task of finding a suitable Lipschitz constant turns out to be rather deceptive for systems in higher dimensions. In contrast, devising a global extended Lipschitz inequality, i.e. (18) is effortless for polynomial vector fields.

In turn, let us recall that behavior of solutions to nonautonomous equation (33) naturally unfolds in 3D-space, $(t, x_1, x_2)$, where the trajectories of this system do not intersect due to uniqueness, but projections of these trajectories on $x_1 - x_2$-plane can intersect with itself and each other. To picture the boundary of stability/trapping region in $(t, x_1, x_2)$-

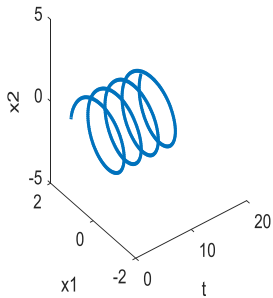
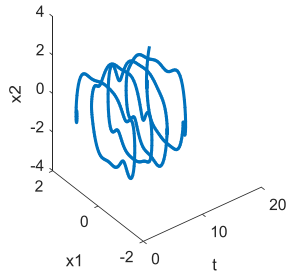
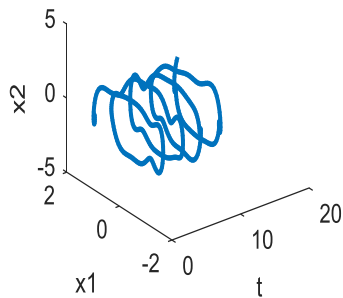

Fig. 3.1    Fig. 3.2    Fig. 3.3

Fig. 3 display 3D-plot of trajectories of system (3.3) computed in reverse time for different values of system's parameters. All trajectories are stemmed from the points locating in interior of stability/trapping regions of (3.3) near their boundaries and rapidly approach the cylinders that are defined the boundaries of the corresponding region in $(t, x_1, x_2)$-space.

space, we plot in Figure 3 the trajectories of system (3.3) computed in *reverse* time for different values of this system's parameters. All trajectories on these figures are stemmed from the points locating in interior of stability/trapping regions of (3.3) and approach the cylinders with axes parallel to $t$-axis that encircle the corresponding regions in $(t, x_1, x_2)$-space. Every other trajectory of system (3.3), simulated in reverse time, approach this cylindrical boundary as well. In Figures 3-6, we set $r_1 = 3.2\pi$, $r_2 = 13$, $\omega_2 = 2\pi$ and $\alpha_2 = -0.05$.

Additionally, in Figure 3.1, $a = 0$, $a_1 = a_2 = 0.1$, whereas in Figure 3.2, $a = 0$, $a_1 = a_2 = 5$ and in Figure 3.3, $a = 0.3$, $a_1 = a_2 = 5$. Clearly, in the first case, the projection of the corresponding trajectory on $x_1 - x_2$-plane yields an oval-like curve, see Figure 4.1a, whereas in the latter two cases, such projections yield fairly irregular bundles, see Figures 4.2a and 4.3a. In each of these cases, the interior regions of these bundles/curve shape the sets of initial vectors, i.e., $x(0) = x_0$, which are included in stability/trapping regions of (33).

Two estimates of these regions, defined by application of formula, $\|W^{-1}(0)x_0\| = X_0$, are shown in dotted-blue and magenta lines in Figures 4.1a – 4.3a. Firstly, $X_0$ is approximated in simulations of the corresponding equation (20). Secondary, $X_0$ is defined as a positive unstable fixed solution of (26) or (31). The second estimate in Fig 4.1a is defined through utility of (26), whereas in Figs. 4.2a and 4.3a such estimate is simulated by utility of the averaged equation (31).

Figure 4.1b – 4.3b plot trajectories of (33) stemmed from initial vectors, which are fitted into the estimated boundaries of stability/trapping regions that are displayed on Figures 4.1a – 4.3a, respectively.

Clearly, simulations of (20) yield a central part of the actual stability basins or trapping regions for this equation.

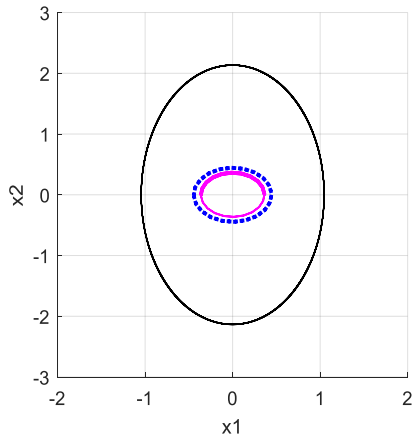

Fig. 4.1a

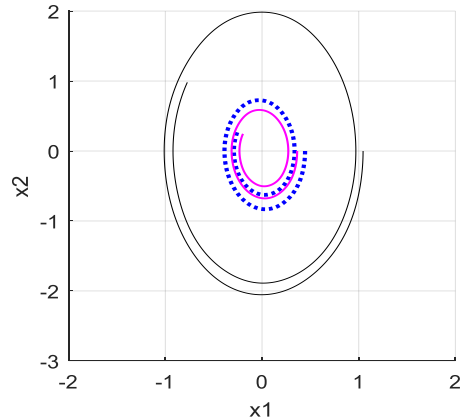

Fig. 4.1b

We notice that the analytical approximations, which define the pertaining initial values of $X_0$ through utility of (26) are close to ones that are obtained in direct simulations of the corresponding equation (20) if $\sup \omega_1(t)$ is relatively small and $a = 0$. If $a \neq 0$, the analytical approximations remain intact if $|a| \ll 1$, since in this case positive roots of (28) sensitively depend upon $a$ and vanish if $a$ increases. Figs. 4.1 show that under these last conditions the pertained values of $x_0$, which are defined either by (26) (blue line) or in direct simulations of (20) (magenta line), are close to each other.

In turn, the averaging technique leads to tolerable analytical estimates of the trapping/stability regions for larger values of both $r_i$ and $a_i$.

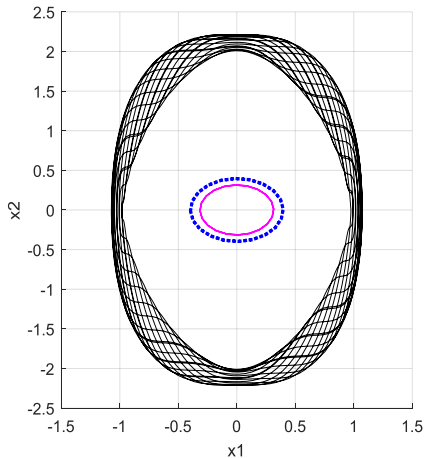

Fig. 4.2a

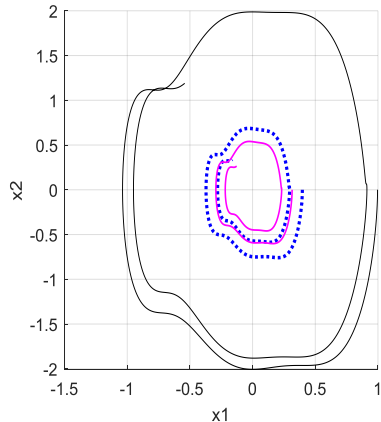

Fig. 4.2b

Figs. 4.2 and 4.3 compare numerical and analytical estimates that are developed in simulations of either equations (20) or (31) corresponding to (33). Clearly, the former two estimates, i.e., magenta and blue lines are sufficiently close to each other for $a_i = 5$ and determine the central part of the actual trapping/stability regions for (33).

Fig.5 plots time-histories of $p(t)$ and $k(t)$, and their running time-averages in blue, yellow, red, and magenta lines, respectively. Both functions notably oscillate, but their running time-averages quickly approach some constant values, which yield the principal contribution to the solutions to (20).

Finally, we apply our methodology for estimating the bounds of solutions to equation (33) with $f = \begin{pmatrix} 0 & -\alpha_2 x_1^3 \end{pmatrix}^T$ and $F = 0$. Such Duffing-like equation is frequently used in simplified modeling of vibrations of elastic structures enforced by oscillations of system's loads [33]. Interestingly, in this case, $\|f\| = \left\|\begin{pmatrix} 0 & -\alpha_2 x_1^3 \end{pmatrix}^T\right\| \leq |\alpha_2| \|x\|^3$ matches the estimate of $\|f\|$ used prior for Van der-Pol-like model. Hence, the auxiliary equations (20) for Duffing-like and Van der-Pol -like models are the same. Fig. 5 plots in blue and red lines time-histories of norms of solutions for Duffing-like and Van der Pol-like models, respectively, emanating from the same initial vector. The black line on this figure plots the estimate of norms of the corresponding solutions, which delivered in simulations of (20). Note that the values of parameters used in Figs. 2 and 5 are identical.

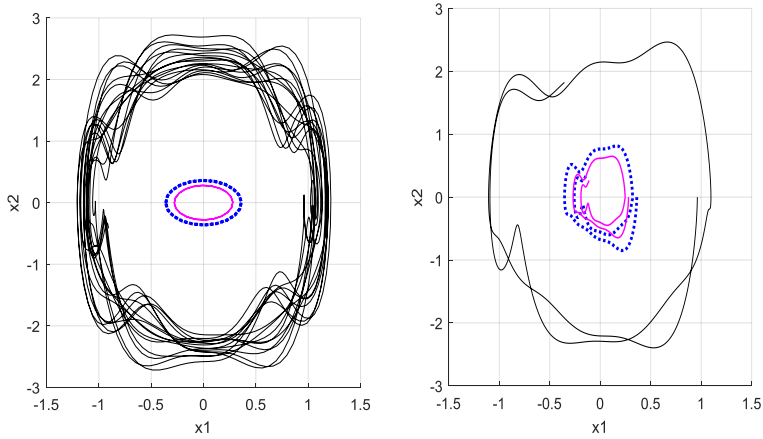

Fig. 4.3a  Fig. 4.3b

Fig. 4. Estimation of stability/trapping regions. Figures 4.1a-4.3a plot approximate boundaries of stability/trapping regions of (33) in projections on $x_1 - x_2$ - plane. Black lines are developed in simulations of (33) in reverse time. Blue-lines are developed in simulations of (20) and magenta-lines are developed through utility of either (26) (Fig. 4.1a) or (31) ( Figs. 4.2a and 4.3a).

Figs 4.1b-4.3b plot projections of trajectories of (33) on $x_1 - x_2$ - plane with initial vectors fitted in the boundary sets that are displayed in Figs. 4.1a-4.3a. On both sets of figures the corresponding lines are matched in color.

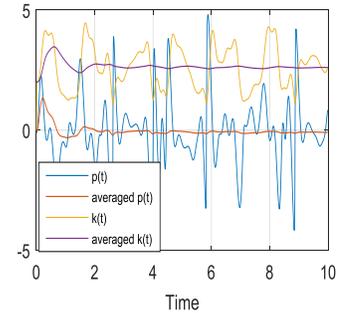

Fig. 5 plot time-histories of $p(t)$ and $k(t)$ in blue and yellow lines and their running time-averages in red and magenta lines, respectively.

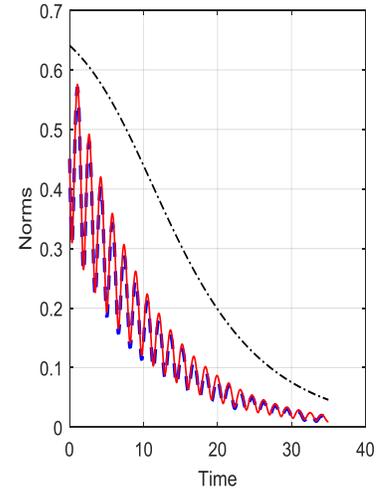

Fig. 6. Estimation of solution norms for Duffing-like and Van der-Pol-like models. Blue, red and black lines plot time-histories of the corresponding norms of solutions and their estimate, respectively.

## 6. Conclusions and Future Work

This paper presents a novel approach to estimation of solution bounds and trapping/stability regions of nonautonomous nonlinear systems. This approach is based on developing of the pivoting differential inequality for the norm of solutions to the initial systems and subsequent analysis of the associated first order auxiliary differential equation. The solutions of this auxiliary equation bound in norm from above the corresponding solutions of the initial systems.

We cast the auxiliary equation in the standard form by using either the Lipschitz condition or its nonlinear extension. The utility of Lipschitz condition linearizes the auxiliary equation and yields the corresponding solution bounds and stability criteria for the conforming nonlinear and nonautonomous systems in the local neighborhoods of the phase space, where the Lipschitz inequality holds. We show that the developed stability criteria turn out to be less conservative than the known ones.

In turn, we developed a nonlinear extension of Lipschitz inequality and applied it to recast the auxiliary equation in a more accurate, but nonlinear and nonautonomous form that, in general, does not admits close form solution. Yet, for autonomous and some other nonlinear systems the solutions of auxiliary equation can be written in close forms.

We formulate the characteristic properties simplifying numerical estimation of the trapping/stability regions of the nonlinear auxiliary equation and consequently apply them for estimation of the corresponding regions for the initial systems. Next, we introduce two approximations reducing the nonlinear auxiliary equation to its autonomous and integrable forms. Analysis of solutions to these autonomous counterparts of the auxiliary equation infer explicit estimates of the trapping/stability regions for the corresponding initial systems, which are contrasted in simulations.

Our theoretical inferences are validated in inclusive numerical simulations that are partly presented in this paper. The simulations show that the accuracy of our estimates inversely correlates with the magnitudes of $\|f(t,x)\|$ and

$\|F(t)\|$, since the auxiliary equation includes only the norms of these perturbations. Hence, the precision of the developed estimates turns out to be adequate if the upper bounds on $\|f\|$ and $\|F\|$ are only known – a frequent premise in theory of systems under uncertainties. But, our approach can yield rather conservative estimates if both $f$ and $F$ are defined precisely.

Yet, the developed approach can be combined with some successive approximations yielding bilateral bounds for the norms of solutions that approach the norms of the accurate solutions under some broad conditions. Application of such refined methodology will be the topic of our subsequent paper.

Conflict of Interest. There is no conflict of interest concerning this paper.

Acknowledgment. This paper was developed in collaboration of its co-authors. The second co-author developed the programs and contributed to interpretation of the simulation results, whereas the first co-author developed the underlined methodology, simulated the pertained systems and drafted the paper.